\documentclass{article}
\usepackage{fullpage}
\usepackage{graphicx}
\usepackage{amsmath}
\usepackage{amssymb}
\usepackage{amsthm}
\usepackage{array}

\newcommand{\T}{\mathbb{R}^3\slash \Lambda}
\newcommand{\Tn}{\mathbb{R}^d \slash \Lambda}

\newcommand{\re}{\operatorname{Re}}

\newtheorem{lemma}{Lemma}[section]
\newtheorem{proposition}[lemma]{Proposition}

\newtheorem{theorem}[lemma]{Theorem}

\newtheorem{question}[lemma]{Question}

\title{The BCC lattice in a long range interaction system}
\author{Xiaofeng Ren
\thanks{Corresponding author.
Phone: 1 202 994-6791; Fax: 1 202 994-6760; E-mail: 
ren@gwu.edu}
\thanks{Supported in part by Simons
Foundation Collaboration Grant for Mathematicians \#709260.}
\\ Department of Mathematics
\\ The George Washington University
\\ Washington, DC 20052, USA
\and Juncheng Wei
\thanks{Supported in part by NSERC-GR008813}
\\ Department of Mathematics
\\ University of British Columbia
\\ Vancouver, BC, Canada  V6T 1Z2
}

\begin{document}
\maketitle

\begin{abstract}
  While the hexagonal lattice is ubiquitous in two dimensions, the body
  centered cubic lattice and the face centered lattice are both commonly
  observed in three dimensions. A geometric variational problem
  motivated by the diblock copolymer theory consists of a short range
  interaction energy and a long range interaction energy. In three dimensions,
  and when the long range interaction is given by the nonlocal operator
  $(-\Delta)^{-3/2}$, it is proved that the body centered cubic lattice is
  the preferred structure. 
\end{abstract}

\section{Introduction}
\setcounter{equation}{0}

In two dimensions the most familiar lattice is the hexagonal lattice,
seen in many places like honeycomb, chicken wire fence, graphene, and carbon
nanotube. In three dimensions, however, there are two common lattices:
the body centered cubic lattice (BCC lattice) and
the face centered cubic lattice (FCC lattice).
In crystallography of metals the BCC lattice is found in iron, chromium,
tungsten, and niobium, while the FCC lattice appears in aluminum, copper,
gold, and silver.
In the sphere packing problem the maximal packing density is attained by 
the hexagonal lattice in two dimensions \cite[T\'oth]{toth},
the FCC lattice in three dimensions \cite[Hales]{hales2},
the $E_8$ lattice in eight dimensions \cite[Viazovska]{viazovska},
and the Leech lattice in 24 dimensions
\cite[Cohn, Kumar, Miller, Radchenko, and Viazovska]
{cohn-kumar-miller-radchenko-viazovska2}. 

Considered in this work is a geometric variational problem whose
free energy functional takes the form
\begin{equation}
  \label{J}
   {\cal J}_{\Lambda,s} (\Omega) = {\cal P}_\Lambda (\Omega) + \frac{\gamma}{2}
   \int_\Omega (-\Delta)^{-s} (\chi_\Omega - \omega) (x) \, dx
\end{equation}
This problem is motivated by 
the Ohta-Kawasaki density functional theory for diblock copolymers \cite{ok}.
A diblock copolymer molecule consists of a
linear subchain of A-monomers grafted covalently to another subchain of
B-monomers \cite[Bates and Fredrickson]{bf}.
Because of the repulsion between the unlike monomers,
the different type sub-chains tend to segregate, but as they are chemically
bonded in chain molecules, segregation of subchains lead to local micro-phase
separation: micro-domains rich in either A-monomers or B-monomers emerge
as a result. The Ohta-Kawasaki theory treats monomer density fields as
the main order parameters. In the
strong segregation regime, the A-monomers occupy a subset $\Omega$ of the
system sample, and the B-monomers occupy the complement of $\Omega$.  
The authors of this paper showed that in this regime the Ohta-Kawasaki
functional converges to a geometric variational problem like \eqref{J}
in the sense of the Gamma limit theory \cite{rw}.

\begin{figure}
\centering
 \includegraphics[scale=0.23]{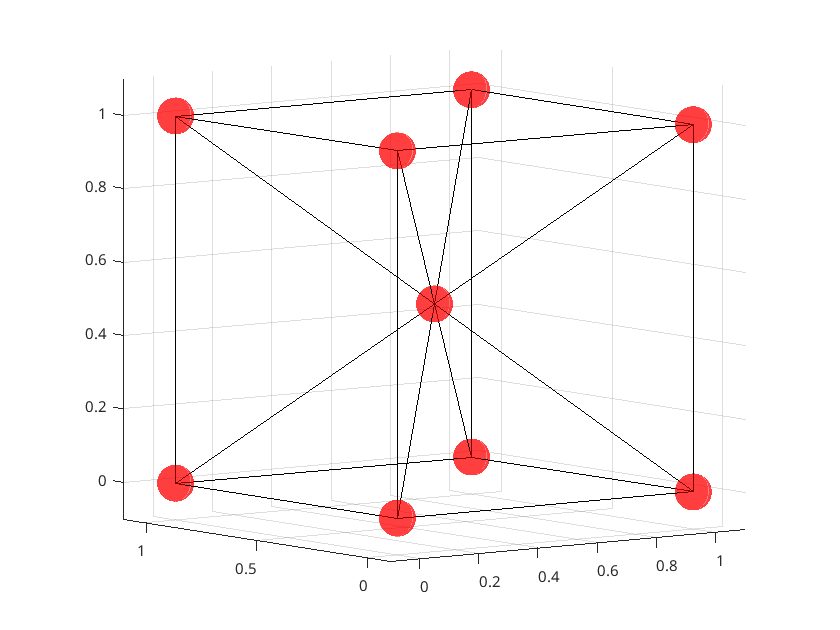} 
 \includegraphics[scale=0.255]{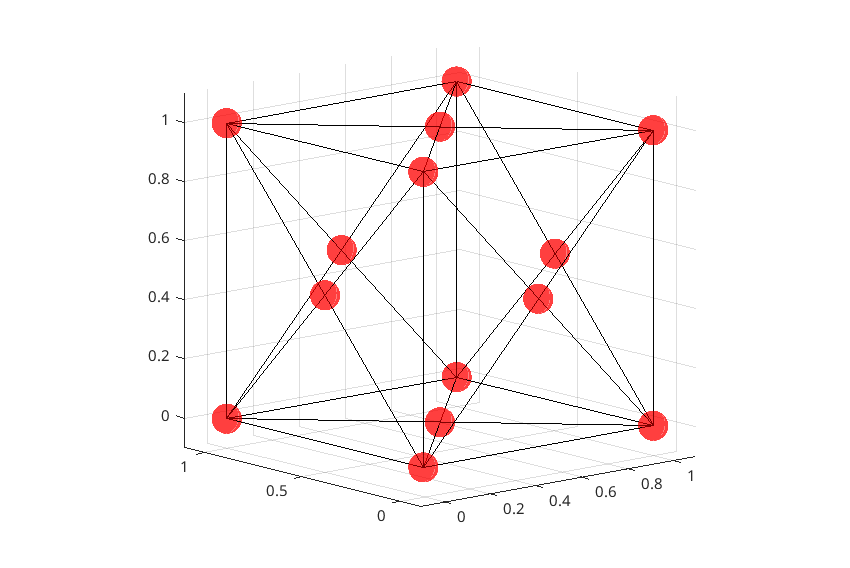}
 \caption{Making a BCC lattice and an FCC lattice out of a cubic lattice.}
\label{f-lattices}
\end{figure}

Generally one cannot expect a perfect periodic structure
if the sample domain has boundary: there are always defects near the boundary.
To circumvent this problem, we take the sample domain to be a flat torus, namely
$\Tn$ where $\Lambda$ is an $d$-dimensional lattice in $\mathbb{R}^d$.
Obviously functions on $\Tn$ (or subsets of $\Tn$) can be viewed as
$\Lambda$-periodic functions on $\mathbb{R}^d$
(or $\Lambda$-periodic subsets of $\mathbb{R}^d$).
Our goal is to find a lattice which is 
optimal in some sense for problem \eqref{J}.

The functional ${\cal J}_{\Lambda,s}$ is defined on Lebesgue measurable
subsets $\Omega$ of $\Tn$ of the prescribed volume:
\begin{equation} \label{omega}
  |\Omega|= \omega |\Tn|, \ \ \mbox{where} \ \ 
  \omega \in (0,1)
\end{equation}
Here $| \cdot |$ denotes the $d$-dimensional Lebesgue measure on $\Tn$.
The number $\omega$ in $(0,1)$ is the first parameter of this problem;
it fixes the volume of $\Omega$. The measure of the flat torus,
denoted $|\mathbb{R}^d/\Lambda |$ or simply $|\Lambda|$,
is the measure of a fundamental parallelepiped of $\Lambda$ in $\mathbb{R}^d$.
It is also called the volume of the lattice $\Lambda$;
see \eqref{absLambda}. 

The first term ${\cal P}_\Lambda(\Omega)$ is the perimeter of $\Omega$
in $\Tn$. If $\Omega$ is a subset with $C^1$ boundary, then the perimeter
is simply the area of $\partial \Omega$. 
The second term in \eqref{J} is an integral over $\Omega$ multipled by
$\gamma/2$ where 
\begin{equation}
  \label{gamma}
  \gamma \in (0,\infty)
\end{equation}
is the second parameter of the
problem. The integrand in the second term is the outcome of the
operator $(-\Delta)^{-s}$ on $\chi_\Omega - \omega$ where
$\chi_\Omega$ is the characteristic function of $\Omega$, i.e.
$\chi_\Omega(x) = 1$ if $x\in \Omega$ and $=0$ if $x \not \in \Omega$. 
The constant $s$
\begin{equation}
  \label{s} s \in (0,\infty)
\end{equation}
in the exponent of $-\Delta$ is the last parameter.

We studied this problem in \cite{rw, rw11, rw12} on bounded domains
with the zero Neumann boundary condition. The exponent $s$ there
was taken to be $1$, so
the operator in the second term was $(-\Delta)^{-1}$ instead of
$(-\Delta)^{-s}$ here.
In the Ohta-Kawasaki theory, $(-\Delta)^{-1}$ is chosen largely for
convenience;  see \cite[Choksi and Ren]{chr}.
In this paper we mainly study the case
$d=3$ and $s=3/2$ and show a connection of this case to
the three dimensional analogy of Kronecker's limit formulas, and to
the height problem for flat tori. 

Denote Green's function of the
operator $(-\Delta)^s$ by $G_{\Lambda,s}$ and decompose $G_{\Lambda,s}$ into
\begin{equation}
  \label{GKR}
   G_{\Lambda,s}(x) = K_{d,s}(x) + R_{\Lambda,s}(x)
\end{equation}
where $K_{d,s}$ is the fundamental solution of the operator
$(-\Delta)^s$ on $\mathbb{R}^d$
and $R_{\Lambda,s}$ is the regular part of $G_{\Lambda,s}$.
While $K_{d,s}$ does not depend on $\Lambda$, 
the regular part of Green's function evaluated at $0$, i.e. $R_{\Lambda,s}I0)$,
contains much information about the lattice $\Lambda$.
To find the optimal lattice, one should minimize $R_{\Lambda, 3/2}(0)$ with
respect to lattice $\Lambda$ of unit volume; see Proposition \ref{p-ItoR}.  

We take the configuration $\Omega$ to be the simplest
set: a ball of radius $r$, i.e. $\Omega =B_r$.
By the constraint \eqref{omega},
$r$ is related to $\omega$ via
\begin{equation}
  \frac{\pi^{d/2} r^d}{\Gamma \left( \frac{d}{2}+1  \right) }
  = \omega |\Tn|
  \end{equation}
Then one compares the free energy per volume,
\begin{equation} \label{e-d}
    \frac{1}{|\Tn|} {\cal J}_{\Lambda, s}(B_r) 
\end{equation}
of the ball on different flat tori $\Tn$. 
It turns out that to minimize \eqref{e-d} it suffices to consider lattices
of unit volume; see Proposition \ref{p-size-shape}. 
Asymptotically, as $r \rightarrow 0$, the quantity \eqref{e-d} is determined by 
$R_{\Lambda, s}(0)$; see Proposition \ref{p-ItoR}. Therefore $R_{\Lambda,s}(0)$
serves as a quantity that measures the
optimality of lattice $\Lambda$ in problem \eqref{J}. 
We raise this question:

\begin{question}
  \label{q-R0}
  Which lattice $\Lambda$ of unit volume minimizes $R_{\Lambda, s}(0)$?
\end{question}

When  $s$ is in $(0,d/2)$,
$K_{d,s}$ in \eqref{GKR} is given by the Riesz potential
\begin{equation}
  K_{d,s}(x) = \frac{\kappa_{d,s}}{|x|^{d-2s}}, \quad 
  \kappa_{d,s} = \frac{\Gamma \left ( \frac{d-2s}{2} \right )}{ 2^{2s} \pi^{d/2}
      \Gamma ( s ) } \label{riesz}
\end{equation}
where $\Gamma(\cdot)$ is the gamma function. When $d=3$ and $s=1$,
the Riesz potential is $\frac{1}{4\pi |x|}$ which is
known as the Coulomb potential.

The borderline case $s=d/2$ turns out to be very interesting, when it comes to 
$R_{\Lambda, d/2}(0)$. If $d=2$, then $s=1$ and $(-\Delta)^s$
is the usual Laplacian $-\Delta$. The fundamental solution $K_{2,1}$ is
not a Riesz potential but a logarithmic function
\begin{equation}
   K_{2,1}(x) = \frac{1}{2\pi} \log \frac{1}{|x|}
\end{equation}
One can write Green's function $G_{\Lambda, 1}$ of $-\Delta$ on a two
dimensional flat torus in terms of a Jacobi's theta function and relate
$R_{\Lambda,1}(0)$ to Dedekind's eta function. 
Chen and Oshita proved that among two dimensional 
lattices of unit area, $R_{\Lambda, 1}(0)$ is minimized
uniquely by the hexagonal lattice \cite{chenoshita-2}.
Sandier and Serfaty noted that one can also connect
$R_{\Lambda,1}(0)$ to Dedekind's eta function
via Kronecker's second limit formula. They gave another proof of
Chen and Oshita's result in \cite{sandier-serfaty}.

In this paper, we consider the case $s=d/2$ in three dimensions, i.e. $d=3$
and $s=3/2$. It will be shown in Lemma \ref{l-1} that the fundamental solution
$K_{3,3/2}$ is
\begin{equation}
   K_{3,3/2}(x) = \frac{1}{2\pi^2} \log \frac{1}{|x|} 
\end{equation}
We will prove the following theorem.

\begin{theorem}
  \label{t-optimal}
  Let $d=3$ and $s=3/2$. Among three dimensional lattices $\Lambda$ of unit
  volume, $R_{\Lambda, 3/2}(0)$ is uniquely minimized by the BCC lattice.
  \end{theorem}

Two lattices are equivalent if we can transform one to the other by
a dilation and an orthogonal transform. Equivalent lattices are
indistinguishable. By a lattice, we often mean the
equivalent class of a lattice. The uniqueness assertion in Theorem
\ref{t-optimal} means uniqueness up to this equivalence. 

The proof of the theorem uses the notion of the height of a manifold.
Let $M$ be a closed Riemannian manifold and  
\begin{equation}
 0=\lambda_0 < \lambda_1 \leq \lambda_2 \leq \lambda_3...
\end{equation}
be the eigenvalues of $-\Delta$ on $M$, 
counting multiplicity.
Define the zeta-regularization
\begin{equation}
 Z(M,s) = \sum_{j=1}^\infty \frac{1}{\lambda_j^s}
\end{equation}
Then the height of $M$ is 
\begin{equation}
  \label{height}
   h(M) = Z'(M,0)
\end{equation}
Here $Z'(M,s)$ is the derivaitve of $Z$ with respect to $s$.

We relate $R_{\Lambda,3/2}(0)$ to $h(\mathbb{R}^3/\Lambda^\ast)$,
the height of $\mathbb{R}^3/\Lambda^\ast$. Here $\Lambda^\ast$ is the dual
lattice of $\Lambda$ (see \eqref{dual} for the definition of a dual lattice),
so $\mathbb{R}^3/\Lambda^\ast$
is again a flat torus viewed as a Riemannian manifold.  We show in
Lemma \ref{l-2} that
\begin{equation}
  \label{h-R}
h(\mathbb{R}^3\slash\Lambda^\ast) = (2\pi)^2 R_{\Lambda,3/2}(0) +2\log(2\pi)
\end{equation}
Then we apply a theorem by Sarnak and Str\"ombergsson:
among lattices in $\mathbb{R}^3$ of unit volume, the height
is minimized uniquely by the FCC lattice \cite{sarnak-strombergsson}.
Since the BCC alttice is the dual lattice of the
FCC lattice, Theorem \ref{t-optimal} follows from \eqref{h-R}. 

\

\noindent {\em Acknowledgments}. We would like to thank Rustum Choksi for
valuable discussions.

\section{Preliminaries}
\setcounter{equation}{0}

The long range interaction term in \eqref{J} is denoted
\begin{equation}
  {\cal I}_{\Lambda, s}(\Omega) = \int_\Omega
  (-\Delta)^{-s} (\chi_\Omega - \omega) (x) \, dx
\end{equation}
so that 
\begin{equation}
  {\cal J}_{\Lambda,s} (\Omega) = {\cal P}_\Lambda (\Omega) 
   + \left ( \frac{\gamma}{2} \right ) {\cal I}_{\Lambda, s}(\Omega)  
\end{equation}
Instead of ${\cal J}_{\Lambda,s}(\Omega)$, it is more appropriate to consider 
\begin{equation}
  \label{tJ}
  \widetilde{\cal J}_{\Lambda,s}(\Omega)
  = \frac{1}{|\Tn|} {\cal J}_{\Lambda,s}(\Omega)
\end{equation}
the energy of the configuration $\Omega$ per volume, when one studies
the impact of the lattice $\Lambda$.

The size and the shape of a configuration play different
roles in $\widetilde{\cal J}_{\Lambda,s}(\Omega)$.
To separate the two factors write the lattice as 
$t \Lambda$ where $ t \in (0, \infty)$ and $\Lambda$ is a lattice of
unit volume, $|\Lambda|=1$.
Then $|t \Lambda|=t^d$. The size of the lattice $t \Lambda$ is given by $t$;
the shape of $t\Lambda$ by $\Lambda$.
The configuration is also written as $t\Omega \subset \mathbb{R}^d / t\Lambda$
with $\Omega \subset \Tn$.

\begin{proposition}
  \label{p-size-shape}
  \begin{enumerate}
  \item For fixed $\Lambda$ and fixed $\Omega$,
    $\widetilde{\cal J}_{t\Lambda,s}(t\Omega)$
    is minimized, with respect to $t$, at
    \[ t=t_{\Lambda,s,\Omega} = \Big ( \frac{{\cal P}_\Lambda (\Omega)}{ s\gamma
      {\cal I}_{\Lambda,s} (\Omega)} \Big )^{1/(2s+1)} \]
    and the minimum value is
    \[  \widetilde{\cal J}_{t_{\Lambda,s,\Omega} \Lambda,s}(t_{\Lambda,s,\Omega}\Omega)
  =    \Big (1+\frac{1}{2s} \Big ) \big (
  {\cal P}_\Lambda (\Omega) \big )^{2s/(2s+1)}
  \big ( s \gamma {\cal I}_{\Lambda,s}(\Omega) \big )^{1/(2s+1)} \]
\item If $\Omega = B_r$ and $\Lambda$ is fixed,
   $\widetilde{\cal J}_{t\Lambda,s}(tB_r)$
   is minimized, with respect to $t$, at
   \[  t=t_{\Lambda,s,B_r} = \Big (\frac{2\pi^{d/2} r^{d-1}}{ s\gamma
      {\cal I}_{\Lambda,s} (B_r) \Gamma(d/2)} \Big )^{1/(2s+1)} \]
   and the minimum value is
   \[ \widetilde{\cal J}_{t_{\Lambda,s,B_r}\Lambda,s}(t_{\Lambda,s,B_r} B_r)
    =    \left (1+\frac{1}{2s} \right )
    \left ( \frac{2\pi^{d/2} r^{d-1}}{\Gamma(d/2)} \right )^{2s/(2s+1)}
        \left ( s \gamma {\cal I}_{\Lambda,s}(B_r) \right )^{1/(2s+1)} \]
  \end{enumerate}
  \end{proposition}

\begin{proof}
We show that
\begin{align}
  {\cal P}_{t\Lambda}(t\Omega) &= t^{d-1} {\cal P}_\Lambda (\Omega)
                                 \label{scaling-local} \\
  {\cal I}_{t\Lambda,s}(t\Omega) &= t^{d+2s} {\cal I}_{\Lambda,
                                   s}(\Omega) \label{scaling-nonlocal}
  \end{align}
The first scaling relation \eqref{scaling-local} follows from the definition
of perimeter:
\begin{equation}
  {\cal P}_\Lambda(\Omega) = \sup \left 
 \{ \int_{\Omega} \mbox{div} \,g(x) \, dx: \ 
 g \in C^1(\Tn, \mathbb{R}^d), \ |g(x)| \leq 1 \ \forall x \in \Tn 
 \right \}.
\end{equation}
To see \eqref{scaling-nonlocal}, let $\lambda_j$, ($j=0,1,2,...$
and $\lambda_0=0< \lambda_1\leq \lambda_2<...$), 
be the eigenvalues of $-\Delta$ on $\Tn$ counting multiplicity, and
$\varphi_j$ be the corresponding eigenfunctions. Let
$v=(-\Delta)^{-s}(\chi_\Omega-\omega)$ and assume
\begin{equation}
  v(x) = \sum_{j=1}^\infty c_j \varphi_j(x), \ x \in \Tn
\end{equation}
Note that $j$ starts from $1$ here. Then 
\begin{equation}
  \chi_\Omega(x) - \omega = (-\Delta)^s v \, (x)
  =  \sum_{j=1}^\infty c_j \lambda_j^s \varphi_j(x)
\end{equation}
On $\mathbb{R}^d / t\Lambda$, the eigenvalues of $-\Delta$ are
$\frac{\lambda_j}{t^2}$, $j=0,1,2,...$, and the corresponding
eigenfunctions are $\varphi_j(\frac{\cdot}{t})$. Then, with
$y \in \mathbb{R}^d / t\Lambda$,
\begin{align*}
 (-\Delta)^s  v\left(\frac{y}{t}\right)
  &= (-\Delta)^s  \sum_{j=1}^\infty c_j \varphi_j\left(\frac{y}{t}\right) \\
  &= \sum_{j=1}^\infty c_j \left( \frac{\lambda_j}{t^2} \right )^s
  \varphi_j\left(\frac{y}{t}\right) \\
 &= t^{-2s} \left ( \chi_{\Omega}\left( \frac{y}{t} \right ) - \omega \right) \\
 &= t^{-2s}  \left ( \chi_{t\Omega} (y) - \omega \right)
\end{align*}
Hence
\begin{equation}
  (-\Delta)^{-s} \left ( \chi_{t\Omega} - \omega \right) (y)
  = t^{2s} v\left(\frac{y}{t}\right)
\end{equation}
and
\begin{equation}
  {\cal I}_{t\Lambda,s}(t \Omega)
  = \int_{t\Omega} t^{2s} v\left(\frac{y}{t}\right) \, dy
  =  t^{n+2s} \int_{\Omega} v(x) \, dx 
  = t^{n+2s} {\cal I}_{\Lambda,s}(\Omega)
  \end{equation} 
which proves \eqref{scaling-nonlocal}.

The energy per cell area of $t \Omega$ is
\begin{align*}
  \nonumber
  \widetilde{\cal J}_{t\Lambda,s}(t\Omega) & =
 \frac{1}{t^n} {\cal J}_{t\Lambda,s}(t \Omega) = 
   \frac{1}{t} {\cal P}_\Lambda (\Omega) 
  +  \frac{t^{2s} \gamma}{2} {\cal I}_{\Lambda,s} (\Omega)
\end{align*}
With respect to $t$, the last quantity is minimized at
\begin{equation}
  \label{t-Omega}
  t=t_{\Lambda,s,\Omega} = \Big ( \frac{{\cal P}_\Lambda (\Omega)}{ s\gamma
    {\cal I}_{\Lambda,s} (\Omega)} \Big )^{1/(2s+1)}
  \end{equation}
and the minimum value is
\begin{equation}
  \label{energey-per-cell}
  \widetilde{\cal J}_{t_{\Lambda,s,\Omega} \Lambda,s}(t_{\Lambda,s,\Omega}\Omega)
  =    \Big (1+\frac{1}{2s} \Big ) \big (
  {\cal P}_\Lambda (\Omega) \big )^{2s/(2s+1)}
        \big ( s \gamma {\cal I}_{\Lambda,s}(\Omega) \big )^{1/(2s+1)}.
\end{equation}
This proves the first part of the proposition.

To see the second part, note that the area of a $(d-1)$-sphere is
\begin{equation}
   {\cal P}_\Lambda (B_r) =\frac{2\pi^{d/2} r^{d-1}}{\Gamma(d/2)}
\end{equation}
which does not depend on $\Lambda$.
\end{proof}

This proposition shows that
to determine the optimal lattice when the configuration set is a ball,
one only needs to minimize ${\cal I}_{\Lambda, s}(B_r)$ with respect to
lattice $\Lambda$ of unit volume. 

The next proposition asserts that as $r \rightarrow 0$,
${\cal I}_{\Lambda, 3/2}(B_r)$ is asymptotically
determined by $R_{\Lambda, 3/2}(0)$. 

\begin{proposition}
  When $d=3$ and $s=3/2$, there exist $c_1(r)$ and $c_2(r)>0$ depending on $r$
    but not on $\Lambda$ such that
    \[ \lim_{r \rightarrow 0} \frac{{\cal I}_{\Lambda, 3/2}(B_r) - c_1(r)}{c_2(r)}
    = R_{\Lambda, 3/2}(0) \]
    for each lattice $\Lambda$
  \label{p-ItoR}
  \end{proposition}

\begin{proof}
  In terms of Green's function
  $G_{\Lambda, s}$ 
\begin{align}
  (-\Delta)^{-s}(\chi_{B_r} - \omega) (w)
  &= \int_{B_r} G_{\Lambda, s} (w-v) \, dv \\
  {\cal I}_{\Lambda,s}(B_r) &= \int_{B_r} \int_{B_r} G_{\Lambda, s} (w-v) \, dv \, dw
  \end{align}

In the case of $d=3$ and $s=3/2$, $G_{\Lambda, 3/2}$ can be written as
\begin{equation}
   G_{\Lambda, 3/2}(w) = -\frac{1}{2\pi^2} \log |w| + R_{\Lambda, 3/2}(w)
\end{equation}
where $R_{\Lambda, 3/2}$ is smooth on
$(\mathbb{R}^3 \backslash \Lambda) \cup \{ (0,0,0) \}$; see Lemma \ref{l-1}.
Note that
\begin{align}
   \int_{B_r} \int_{B_r} \log |w-\tilde{w}| \, dw d\tilde{w} &=
   \left( \frac{4\pi}{3} \right)^2 r^6 \log \frac{1}{r}
   - r^6 \int_{B_1} \int_{B_1} \log |v -\tilde{v}| \, dv d \tilde{v} \\
   \int_{B_r} \int_{B_r} R_{\Lambda, 3/2}(w-\tilde{w})  \, dw d\tilde{w} &=
   \int_{B_r} \int_{B_r} \left ( R_{\Lambda, 3/2}(0) +O(r^2) \right )
   \, dw d\tilde{w} \nonumber \\
   &= \left( \frac{4\pi}{3} \right)^2 r^6 R_{\Lambda, 3/2}(0) + O(r^8)
\end{align}
Here one used the fact $\nabla R(0)=0$ since $G(w) = G(-w)$ for all
$w \in \mathbb{R}^3 \backslash \Lambda$ to deduce the $O(r^2)$ term.
Therefore
\begin{equation}
  \label{IR0}
  {\cal I}_{\Lambda, 3/2}(B_r) =
   \frac{8 r^6}{9} \log \frac{1}{r} - \frac{r^6}{2\pi^2}
  \int_{B_1} \int_{B_1} \log |v -\tilde{v}| \, dv d \tilde{v}
   + \left( \frac{4\pi}{3} \right)^2 r^6 R_{\Lambda,3/2}(0) + O(r^8)
\end{equation}
Let
\begin{align*} c_1(r) &= \frac{8 r^6}{9} \log \frac{1}{r} - \frac{r^6}{2\pi^2}
\int_{B_1} \int_{B_1} \log |v -\tilde{v}| \, dv d \tilde{v} \\
c_2(r) &= \left( \frac{4\pi}{3} \right)^2 r^6
\end{align*}
and the proposition follows from \eqref{IR0}.
\end{proof}

In the rest of the paper we show that $R_{\Lambda, 3/2}(0)$ is uniquely
minimized by the BCC lattice.  Let us briefly recall some basic
facts in the lattice theory.
A $d$-dimensional  lattice $\Lambda$ is a subset of $\mathbb{R}^d$ of the form 
\begin{equation}
 \Lambda = \{ n_1 v_1+ n_2 c_2+...+ n_d v_d: \ n_j \in \mathbb{Z} \}
\end{equation}
where $\{ v_1, v_2, ..., v_d\}$ is a set of linearly independent vectors
in $\mathbb{R}^d$. Denote by $V \in GL(d,\mathbb{R})$
the $d \times d$ matrix whose $j$-th row
is the row vector $v_j$; $V$ is called a generator matrix of $\Lambda$.
The Gram matrix $Q$ of $V$ is a positive definite matrix given by
\begin{equation}
  Q=V V' 
\end{equation}
where $V'$ is the transpose of $V$.

The most important lattice in two dimensions is arguably the hexagonal lattice.
For this lattice we can take $v_1=(1,0)$ and $v_2=(1/2, \sqrt{3}/2)$. Then 
\begin{equation}
  \label{h-lattice}
  V= \left [  \begin{array}{cc} 1 & 0 \\ \frac{1}{2} & \frac{\sqrt{3}}{2}
    \end{array} \right ], \
  Q=  \left [ \begin{array}{cc} 1 & \frac{1}{2} \\ \frac{1}{2} & 1
    \end{array} \right ] 
\end{equation}

When $d=3$, $\mathbb{Z}^3$ is a cubic lattice, generated by $(1,0,0)$,
$(0,1,0)$, and $(0,0,1)$. Both the generator matrix and the Gram matrix are
both the identity matrix. 

Adding {\em body centers}
$(1/2,1/2,1/2) + \lambda$, $\lambda \in \mathbb{Z}^3$, to $\mathbb{Z}^3$,
we obtain a BCC lattice; see the left plot of Figure \ref{f-lattices}.
This lattice has a
generator matrix and a corresponding Gram matrix as follows:
\begin{equation}
  \label{bcc-lattice}
  V= \left [  \begin{array}{ccc} -\frac{1}{2} & \frac{1}{2} & \frac{1}{2}
      \\ \frac{1}{2} & -\frac{1}{2} & \frac{1}{2}
      \\ \frac{1}{2} & \frac{1}{2} & -\frac{1}{2} \end{array} \right ], \
  Q=  \left [ \begin{array}{ccc} \frac{3}{4} & -\frac{1}{4} & -\frac{1}{4}
      \\ -\frac{1}{4} & \frac{3}{4}  & -\frac{1}{4}
      \\ -\frac{1}{4} & -\frac{1}{4} & \frac{3}{4} \end{array} \right ]
\end{equation}

From the cubic lattice $\mathbb{Z}^3$, adding {\em face centers}
$(0,1/2,1/2) + \lambda$, $(1/2,0,1/2)+\lambda$, and $(1/2,1/2,0)+\lambda$,
$\lambda \in \mathbb{Z}^3$, we have an FCC lattice; see the right plot
of Figure \ref{f-lattices}.  This lattice has
\begin{equation}
  \label{fcc-lattice}
  V= \left [  \begin{array}{ccc} 0 & \frac{1}{2} & \frac{1}{2}
      \\ \frac{1}{2} & 0 & \frac{1}{2} \\ \frac{1}{2} & \frac{1}{2} & 0
    \end{array} \right ], \
  Q=  \left [ \begin{array}{ccc} \frac{1}{2} & \frac{1}{4} & \frac{1}{4}
      \\ \frac{1}{4} & \frac{1}{2}  & \frac{1}{4}
      \\ \frac{1}{4} & \frac{1}{4} & \frac{1}{2} \end{array} \right ]
\end{equation}
as a generator matrix and a corresponding Gram matrix.

The flat torus associated with $\Lambda$ is the quotient space $\Tn$. 
The Lebesgue measure of $\Tn$ is 
\begin{equation}
  \label{absLambda}
  |\Tn|=|\Lambda| = |\det V| = (\det Q)^{1/2}
\end{equation}
which we denote by $|\Lambda|$ for simplicity. 
It is also called the volume of the lattice $\Lambda$.

Two lattices are equivalent if one can be transformed to the other by
a dilation and an orthogonal transform.
Any lattice that is equivalent to the lattice described in \eqref{h-lattice},
(\eqref{bcc-lattice} or \eqref{fcc-lattice}, respectively),
is a hexagonal lattice, (BCC lattice or FCC lattice, respectively). 

To understand this equivalence relation in terms of generator matrices, let 
$V \in GL(d,\mathbb{R})$ transform to  $\tilde{V} \in GL(d,\mathbb{R})$ by
the right action
\begin{equation}
 V \rightarrow    \tilde{V} = V \kappa U 
\end{equation}
where $\kappa \in \mathbb{R}^\times = \mathbb{R} \backslash \{0\}$
and $U$ is a $d\times d$ orthogonal matrix and.
The resulting space of left cosets is
\begin{equation}
  \mathbb{H}^d = GL(d,\mathbb{R}) / O(d) \mathbb{R}^\times
\end{equation}
where $O(d)$ is the group of $d\times d$ orthogonal matrices.

When $d=3$, by the Iwasawa decomposition \cite{iwasawa},
each point in $\mathbb{H}^3$ can
be uniquely represented by an upper triangular matrix of the form
\begin{equation}
\tau = \left (  \begin{array}{ccc} y_1y_2 & y_1x_2 & x_3
    \\ & y_1 & x_1 \\ & & 1 \end{array} \right ) \label{tau}
\end{equation}
where $x_1,x_2,x_3,y_1,y_2 \in \mathbb{R}$ and $y_1, y_2 >0$. More precisely,
every generator
matrix $V \in GL(3, \mathbb{R})$ can be written as
\begin{equation}
  V= \tau \kappa U \label{Vtau}
\end{equation}
where $\tau$ is of the form \eqref{tau}, $\kappa \in \mathbb{R}^\times$, and
$U$ is an orthogonal matrix. The Gram matrix $Q$ corresponding to $V$ is
\begin{equation}
  Q = \kappa^2 \tau \tau' \label{Qtau}
  \end{equation}

The generator matrix $V\in GL(3,\mathbb{R})$ of the BCC lattice in
\eqref{bcc-lattice} is decomposed as in \eqref{Vtau} with 
\begin{equation}
  \tau =\left [\begin{array}{ccc}
      \frac{2\sqrt{2}}{3} \ \frac{\sqrt{3}}{2} & \frac{2\sqrt{2}}{3} \
      (-\frac{1}{2}) & -\frac{1}{3} \\
      & \frac{2\sqrt{2}}{3}  & -\frac{1}{3} \\
      &   &  1 \end{array} \right ]
    , \ \kappa=\frac{\sqrt{3}}{2}, \ U=\left [\begin{array}{ccc}
    0 & \frac{1}{\sqrt{2}} & \frac{1}{\sqrt{2}} \\
    \frac{2}{\sqrt{6}} & -\frac{1}{\sqrt{6}} & \frac{1}{\sqrt{6}} \\
    \frac{1}{\sqrt{3}} & \frac{1}{\sqrt{3}} & -\frac{1}{\sqrt{3}}
    \end{array} \right ]
\end{equation}
For the FCC lattice in \eqref{fcc-lattice} the decomposition \eqref{Vtau}
is given by
\begin{equation}
  \ \tau =\left [\begin{array}{ccc}
      \frac{\sqrt{3}}{2} \ \frac{2\sqrt{2}}{3} & \frac{\sqrt{3}}{2} \
      \frac{1}{3} & \frac{1}{2} \\
     & \frac{\sqrt{3}}{2}  & \frac{1}{2} \\
     &  &  1 \end{array} \right ]
    , \ \kappa=\frac{1}{\sqrt{2}}, \ U=\left [\begin{array}{ccc}
    -\frac{1}{\sqrt{3}} & \frac{1}{\sqrt{3}} & \frac{1}{\sqrt{3}} \\
    \frac{1}{\sqrt{6}} & -\frac{1}{\sqrt{6}} & \frac{2}{\sqrt{6}} \\
    \frac{1}{\sqrt{2}} & \frac{1}{\sqrt{2}} & 0
    \end{array} \right ]
\end{equation}

There is still more redundancy in $\mathbb{H}^d$. 
If $g$ is in the modular group
$SL(d,\mathbb{Z})$ and $V$ is a generator matrix, then
$g V$ is another generator matrix of the same lattice.
This left action, $V \rightarrow gV$, gives rise to the
double coset space
\begin{equation}
  SL(d,\mathbb{Z}) \backslash \mathbb{H}^d
  =  SL(d,\mathbb{Z}) \backslash   GL(d,\mathbb{R}) / O(d)  \mathbb{R}^\times
\end{equation}
which is the space of equivalent lattice classes. 

If $\Lambda$ is a lattice in $\mathbb{R}^d$,  the dual lattice of $\Lambda$
is 
\begin{equation}
  \label{dual}
  \Lambda^\ast = \{ k \in \mathbb{R}^d: k \cdot \lambda \in \mathbb{Z}
   \ \mbox{for all} \ \lambda \in \Lambda \}
\end{equation}
If $V$ and $Q$ are a geneorator and a
Gram matrices of $\Lambda$ respectively, then we can take
\begin{equation}
   V^\ast = (V^{-1})', \ \ Q^\ast = V^\ast (V^{\ast})' = Q^{-1}
\end{equation}
to be a generator and a Gram matrices of $\Lambda^\ast$ respectively.
The dual lattice of a hexagonal lattice is again a hexagonal lattice.
The dual lattice of an FCC lattice is a BCC lattice, and vice versa.

As a function of $w \in \mathbb{R}^d$,
$e^{2\pi i k\cdot w}$ is $\Lambda$-periodic if
$k \in \Lambda^\ast$; it is an eigenfunction of the $-\Delta$ operator on
$\Tn$, with the eigenvalue $(2\pi)^2 |k|^2$. For $s>0$
let $G_{\Lambda,s}$ be Green's function of the operator $(-\Delta)^s$ on $\Tn$ so that
\begin{equation}
  \label{G}
  (-\Delta)^s G_{\Lambda,s} = \delta - \frac{1}{|\Lambda|},
  \  \int_{\Tn} G_{\Lambda,s}(w) \, dw =0
  \end{equation}
Now write $G_{\Lambda,s}$ as
\[ G_{\Lambda,s}(w) = \sum_{k \in \Lambda^\ast} C_k e^{2\pi i k\cdot w},
\ C_0 =0\]
where the sum is assumed to converge in the sense of distributions on $\Tn$.
Here $C_0=0$ because of the integral condition in \eqref{G}.
By the fact
\begin{equation}
  \delta = \sum_{k \in \Lambda^\ast} \frac{1}{|\Lambda|} e^{2\pi i k \cdot w},
\end{equation}
the $C_k$'s satisfy
\[ \sum_{k \in \Lambda^\ast} (2\pi)^{2s} |k|^{2s} C_k e^{2\pi i k\cdot w}
= \sum_{k \in \Lambda^\ast}  \frac{1}{|\Lambda|} e^{2\pi i k\cdot w}
 - \frac{1}{|\Lambda|}.
\]
Hence $C_k = \frac{1}{(2\pi)^{2s} |k|^{2s} |\Lambda|}$
if $k \ne 0$ and $C_0=0$. Consequently
\begin{equation}
  G_{\Lambda,s}(w) = \sum_{k \in \Lambda^\ast\backslash \{0 \}}
  \frac{ e^{2\pi i k\cdot w}}
  {(2\pi)^{2s} |k|^{2s} |\Lambda| } \label{Green-Fourier}
\end{equation}

Let $V^\ast$ be a generator matrix of $\Lambda^\ast$ and $Q^\ast =
V^\ast (V^\ast)'$ be the associated Gram matrix. If $k= m V^\ast$ for
$m \in \mathbb{Z}^d$, then 
\[ |k|^2 = m V^\ast ( m V^\ast )' = m Q^\ast m' := Q^\ast [m] \]
This allows us to write alternatively
\begin{equation}
  \label{G2}
  G_{\Lambda,s}(w)= \sum_{m \in \mathbb{Z}^d \backslash \{ 0 \}}
  \frac{e^{2\pi i m  V^\ast \cdot w}}{(2\pi)^{2s} Q^\ast[m]^s |\Lambda|}
  \end{equation}

\section{$R_{\Lambda, 3/2}(0)$}
\setcounter{equation}{0}

The expression \eqref{G2} for Green's function $G_{\Lambda,s}$ is reminiscent
of the twisted Epstein zeta function introduced by
Siegel. Let $Q$ be a positive definite $d\times d$ matrix,
$u, v \in \mathbb{R}^d$, and define this zeta function by 
\begin{equation}
  \zeta(s,u,v,Q) = \sum_{m \in \mathbb{Z}^d, m+v \ne 0}
  \frac{e^{2\pi i m\cdot u}}{ Q[m+v]^s}, \  \
    s\in \mathbb{C}, \ \re s >\frac{d}{2}.
  \label{zeta}
\end{equation}

\begin{lemma}[{\cite[Epstein]{epstein}}, {\cite[Siegel]{siegel}}]
  As a function of $s$, $\zeta$ admits an analytic continuation on
  $\mathbb{C}$. It is entire if $u \not \in \mathbb{Z}^d$,
  and has a simple pole at $s=d/2$ if $u \in \mathbb{Z}^d$. 
Furthermore, $\zeta$ satisfies the functional equation
\begin{equation}
  \pi^{-s} \Gamma(s) \zeta(s,u,v,Q) = e^{-2\pi i u \cdot v}
  (\det Q)^{-1/2} \pi^{-(d/2-s)} \Gamma\left (\frac{d}{2} -s \right)
  \zeta\left(\frac{d}{2}-s, v, -u, Q^{-1} \right ).
\end{equation}
\end{lemma}

The case $u \in \mathbb{Z}^d$ was treated by Epstein and the case
$u \not \in \mathbb{Z}^d$ by Siegel. In the latter case there is 
the question what $\zeta(d/2, u, v, Q)$ is? If $d=2$, the
answer is given by Kronecker's second limit formula. 
If $d=3$ and $v=0$, Efrat proved the result below.

By \eqref{Vtau} and \eqref{Qtau} we write
\begin{equation}
  \label{Qtau2}
  Q  = (\det Q)^{1/3} (y_1^2y_2)^{-2/3} \tau \tau'
\end{equation}
where $\tau$, together with the $x_j$'s and the $y_j$'s,
is given in \eqref{tau}. Let 
\begin{equation}\label{z}
  z_1 = x_1+i y_1, \ z_2=x_2 + i y_2, \ x_4 = x_3-x_1x_2. \end{equation}
Define a matrix $Q_2$ to be
\begin{equation}
Q_2= \frac{1}{y_1} \left ( \begin{array}{cc} y_1 & x_1 \\ & 1 \end{array}
\right ) \left (  \begin{array}{cc} y_1 & \\ x_1 & 1 \end{array} \right ) 
\end{equation}
and denote by $\zeta_2$ this zeta function in two dimensions:
\begin{align}
  \zeta_2 \left (s, \left [ \begin{array}{c} u_2 \\ u_3 \end{array} \right ],
  0, Q_2 \right )
  & = \sum_{(n_1,n_2) \in \mathbb{Z}^2 \backslash \{(0,0)\}} \frac{e^{2\pi i
    (n_1 u_2+n_2u_3)}}{Q_2[n]^s} \nonumber
 \\ & = \sum_{(n_1,n_2) \in \mathbb{Z}^2 \backslash \{(0,0)\}}
    e^{2\pi i (n_1 u_2+n_2 u_3)}
      \frac{y_1^s}{|n_1z_1+n_2|^{2s}}
\end{align}
which converges absolutely for $\re s >1$.

\begin{lemma}[{\cite[Efrat]{efrat}}]
    Let $d=3$, $v=0$, and $ u \not \in \mathbb{Z}^3$. Then
\label{l-kronecker-2}
  \begin{align*}
  (\det Q)^{1/2} \zeta \left( \frac{3}{2}, u, 0, Q \right )
 & = y_1^{1/2} y_2 \ \zeta_2\left ( \frac{3}{2}, \left [\begin{array}{l} u_2 \\ u_3
    \end{array}\right ], 0 , Q_2 \right )  \\
  & \quad  - 4\pi \log \prod_{(n_1,n_2) \in \mathbb{Z}^2}
   \Big | 1- \exp \Big (-2\pi y_2 |(n_2-u_3)z_1 - (n_1-u_2)| \Big )  \\
   & \quad \exp \Big ( 2\pi i (u_1 + (n_1-u_2)x_2+ (n_2-u_3)x_4) \Big )  \Big |  
\end{align*}
\end{lemma}

A comparison of $G_{\Lambda,s}$ in \eqref{G2} and $\zeta$ in \eqref{zeta}
shows that
\begin{equation}
  \label{Gzeta}
  G_{\Lambda,s}(w)
  = (2\pi)^{-2s} |\Lambda|^{-1} \zeta \left (s, V^\ast w, 0, Q^\ast \right )
  \end{equation}
We proceed to find the fundamental solution of $G_{\Lambda, 3/2}$ and
the regular part of  $G_{\Lambda, 3/2}$ evaluated at $0$. Let $V^\ast$
be a generator matrix of the dual lattice $\Lambda^\ast$ and
$Q^\ast=V^\ast (V^\ast)'$ be its Gram matrix. By the Iwasawa decomposition
one can write $V^\ast$ and $Q^\ast$ as 
\begin{align}
  V^\ast&= (\det V^\ast)^{1/3} (y_1^2 y_2)^{-1/3} \tau U \label{iwasawa-Vast}\\
  Q^\ast &= (\det Q^\ast)^{1/3} (y_1^2y_2)^{-2/3} \tau \tau' \label{iwasawa-Q-ast} 
  \end{align}
Here $U$ is an orthogonal matrix and
$\tau$ is a triangular matrix in the Iwasawa decomposition of $V^\ast$
(not $V$); $\tau$ takes the form 
\begin{equation}
  \label{Vast-tau}
\tau = \left (  \begin{array}{ccc} y_1y_2 & y_1x_2 & x_3
    \\ & y_1 & x_1 \\ & & 1 \end{array} \right )
\end{equation}
where $y_1,y_2>0$, and $x_1,x_2,x_3 \in \mathbb{R}$. Also
set $z_1$, $z_2$, and $x_4$ as in \eqref{z}.

\begin{lemma} \label{l-1}
  When $d=3$ and $s = 3/2$, 
  \[ G_{\Lambda,3/2}(w) = - \frac{1}{2 \pi^2} \log |w| + R_{\Lambda,3/2}(w) \]
  where $R_{\Lambda,3/2}$ is a smooth function near $0$ and
  \begin{align*}
    R_{\Lambda,3/2}(0) &= \frac{y_1^{1/2} y_2}{(2\pi)^{3}}
    \left (  \sum_{(n_1,n_2) \in \mathbb{Z}^2\backslash \{(0,0)\} }
  \frac{y_1^{3/2}}{|n_1z_1+n_2|^{3}} \right )
    - \frac{1}{2\pi^2}  \log \left ( 2\pi y_1^{1/3} y_2^{2/3} |\Lambda|^{-1/3} \right) \\
    & \quad 
    -\frac{1}{2\pi^2} \log \prod_{(n_1,n_2) \in \mathbb{Z}^2 \backslash \{(0,0)\}}
    \Big | 1- \exp \Big (-2\pi y_2 |n_2 z_1 - n_1| 
    + 2\pi i ( n_1 x_2+ n_2 x_4) \Big )  \Big |
  \end{align*}
  Here $x_j$, $y_j$, and $z_j$, given by \eqref{Vast-tau} and
  \eqref{z}, come from $\tau$ in the Iwasawa decomposition of $V^\ast$.
 \end{lemma}

\begin{proof}
By \eqref{Gzeta} and Lemma \ref{l-kronecker-2}
\begin{align}
  G_{\Lambda,3/2}(w) &= (2\pi)^{-3} (\det Q^\ast)^{1/2}
  \zeta \left (3/2, V^\ast w, 0, Q^\ast \right ) \\
  &= \frac{y_1^{1/2} y_2}{(2\pi)^3} \ \zeta_2
    \left ( \frac{3}{2}, \left [\begin{array}{l} u_2 \\ u_3
    \end{array}\right ], 0 , Q_2^\ast \right )  \nonumber \\
  & \quad  - \frac{1}{2\pi^2} \log \prod_{(n_1,n_2) \in \mathbb{Z}^2}
   \Big | 1- \exp \Big (-2\pi y_2 |(n_2-u_3)z_1 - (n_1-u_2)| \Big )\nonumber \\
   & \quad \exp \Big ( 2\pi i (u_1 + (n_1-u_2)x_2+ (n_2-u_3)x_4) \Big )  \Big |
   \label{ttt}
\end{align}
Here
\begin{equation}
  \label{uw}
  u= V^\ast w 
\end{equation}
and
\begin{equation}
Q_2^\ast = \frac{1}{y_1} \left ( \begin{array}{cc} y_1 & x_1 \\ & 1 \end{array}
\right ) \left (  \begin{array}{cc} y_1 & \\ x_1 & 1 \end{array} \right ) 
\end{equation}

As $w \rightarrow 0$ in \eqref{ttt}, $u \rightarrow 0$; since  
the function $\zeta_2$ is regular at $u_2=u_3=0$,
\begin{align}
 \zeta_2\left (\frac{3}{2}, 0, 0, Q^\ast_2 \right) &= \sum_{(n_1,n_2)\in \mathbb{Z}^2\backslash\{(0,0)\}} \frac{y_1^{3/2}}{|n_1z_1+n_2|^3}  \label{zeta2-term}
\end{align}
The term on the right side of \eqref{zeta2-term}
is actually the real analytic Eisenstein series at $s=3/2$ times twice the
Riemann zeta function at $3$.

In the infinite product of \eqref{ttt} the
$n_1=n_2=0$ term 
\begin{align}
  \Big |  1- \exp \Big (-2 \pi y_2 | - u_3 z_1 + u_2 | + 2\pi i (u_1 - u_2 x_2 - u_3 x_4)  \Big ) \Big | \label{singular} 
\end{align}
causes a logarithmic singularity at $w=0$. 
Let
\begin{equation}
  \tilde{w} = U w
\end{equation}
where $U$ is the orthogonal matrix in \eqref{iwasawa-Vast}. Then
by \eqref{iwasawa-Vast} and \eqref{uw}
\begin{align} u &= (\det V^\ast)^{1/3} (y_1^2y_2)^{-1/3} \tau U w
  \nonumber \\
   &= (\det V^\ast)^{1/3} (y_1^2y_2)^{-1/3} \tau \tilde{w}
\end{align}
namely
\begin{align*}
  u_1 &=(\det V^\ast)^{1/3}  (y_1^2y_2)^{-1/3}
  (y_1 y_2 \tilde{w}_1 + y_1 x_2 \tilde{w}_2 + x_3 \tilde{w}_3) \\
  u_2 &=(\det V^\ast)^{1/3} (y_1^2y_2)^{-1/3} (y_1 \tilde{w}_2 + x_1 \tilde{w}_3)
   \\
  u_3 &=(\det V^\ast)^{1/3}  (y_1^2y_2)^{-1/3} \tilde{w}_3
\end{align*}
Calculations show that
\begin{align}
  q &:=-2 \pi y_2 | - u_3 z_1 + u_2 | + 2\pi i (u_1 - u_2 x_2 - u_3 x_4) \\
  & = 2\pi y_1^{1/3} y_2^{2/3} \left ( -|\tilde{w}_2- i \tilde{w}_3|
  + i \tilde{w}_1 \right ) (\det V^\ast)^{1/3}  \label{tt}
\end{align}
Then \eqref{tt} implies that
\begin{align}
  \log | 1- e^q| &= \log |q| + O(|q|) \nonumber \\
  &= \log |\tilde{w}| + \log \left ( 2\pi y_1^{1/3} y_2^{2/3} \right)  
  + \log |\det V^\ast|^{1/3} 
  + O \left (y_1^{1/3}y_2^{2/3} |\tilde{w}|\ |\det V^\ast|^{1/3} \right )
  \nonumber \\
  &=\log |w| + \log \left ( 2\pi y_1^{1/3} y_2^{2/3} \right)  
  + \log |\det V^\ast|^{1/3} 
  + O \left (y_1^{1/3}y_2^{2/3} |w|\ |\det V^\ast|^{1/3} \right )
  \label{log-term}
  \end{align}

The remaining terms in the infinite product of \eqref{ttt}
are regular; they give
\begin{align}
  \prod_{(n_1,n_2) \in \mathbb{Z}^2\backslash \{(0,0)\}}
  \Big | 1- \exp \Big (-2\pi y_2 |n_2 z_1 - n_1|
   + 2\pi i (n_1 x_2+ n_2 x_4) \Big )  \Big |  \label{product-term}
  \end{align}
when $u=0$, i.e. $w=0$. Lemma \ref{l-1} then follows from
\eqref{ttt} with the help of
\eqref{zeta2-term}, \eqref{log-term}, and \eqref{product-term}.
\end{proof}

\section{$k_3(Q^\ast)$ and $h(\mathbb{R}^3/\Lambda^\ast)$}
\setcounter{equation}{0}

Now consider Epstein's zeta function
\begin{equation}
  Z_Q(s) = \sum_{m \in \mathbb{Z}^d\backslash \{0\}} \frac{1}{Q[m]^s}, \
  s \in \mathbb{C}, \ \ \re s > \frac{d}{2}
\end{equation}
where $Q$ is a positive definite $d \times d$ matrix.
Clearly $Z_Q$ is a special case of the twisted Epstein zeta function considered
earlier:
\begin{equation}
  Z_Q(s) = \zeta(s, u, 0, Q), \ \mbox{if} \ u \in \mathbb{Z}^d
\end{equation}

\begin{lemma}[{\cite[Epstein]{epstein}}]
\label{l-functional}
  As a function of $s$, $Z_Q(s)$  has a meromorphic extension to
  $\mathbb{C}$ with only one
  simple pole at $d/2$ and its residue is
  $(\det Q)^{-1/2} \pi^{d/2} \Gamma(d/2)^{-1}$.
  Moreover, $Z_Q$ satisfies the functional equation
    \[ \pi^{-s} \Gamma(s)Z_Q (s) = (\det Q)^{-1/2} \pi^{-(d/2-s)}
    \Gamma\left (\frac{d}{2}-s \right) Z_{Q^{-1}} \left(\frac{d}{2}-s \right) \]
  \end{lemma}

 Let $k_d(Q)$ be
 \begin{equation}
   \label{k}
   k_d(Q)= \lim_{s\rightarrow d/2} \left ( Z_Q(s) -
   \frac{(\det Q)^{-1/2} \pi^{d/2} \Gamma(d/2)^{-1}}{s-d/2} \right )
 \end{equation}
the constant term in the Laurent series of $Z_Q$ about $d/2$. 
When $d=2$, $k_2(Q)$ is given by  the first Kronecker limit formula.
When $d=3$, $k_3(Q)$ follows from a formula in \cite[Efrat]{efrat2}.
    
\begin{lemma} \label{l-k3}
  Let $d=3$ and $\tau$ be the triangular matrix associate with $Q$
  as in \eqref{tau}, \eqref{Vtau} and \eqref{Qtau}. Then
     \[ k_3(Q)= 2\pi (\det Q)^{-1/2} \left ( \log(\det Q)^{-1/3} + 2 \gamma -2 -
       2 \log y_1^{1/3} y_2^{2/3}
       - 4 \log g(\tau)  \right ) \]
     where $\gamma=0.57721...$ is Euler's constant and  
     \begin{align}
       g(\tau) &= \exp\left (-\frac{y_1^{1/2}y_2}{8\pi}
       \left ( \sum_{(n_1,n_2) \in \mathbb{Z}^2\backslash\{(0,0)\}}
       \frac{y_1^{3/2}}{|n_1 z_1 + n_2|^3 } \right )
       \right ) \nonumber \\ 
       & \quad \prod_{(n_1,n_2) \in (\mathbb{Z}^2\backslash
         \{(0,0)\})/\pm 1} \Big | 1- \exp \Big (-2\pi y_2 |n_2 z_1 - n_1| 
   + 2\pi i ( n_1 x_2+ n_2 x_4) \Big )  \Big | \label{gtau}
       \end{align}
  \end{lemma}

\begin{proof}

As explained in \cite{efrat2}, the Eisenstein series $E(\tau,t)$
associated with the maximal parabolic subgroup of $SL(3,\mathbb{Z})$
can be written as
\begin{equation}
    \label{E}
    E(\tau,t)= \sum_{(n_1,n_2,n_3)=1} \frac{(y_1^2y_2)^t}{
    \left ( y_1^2|n_3 z_2+n_1|^2 + (n_3 x_3 + n_1 x_1 +n_2)^2\right )^{3t/2}}
\end{equation}
where $x_j$, $y_j$, and $z_j$ come from $\tau$ as in \eqref{tau}
and \eqref{z}. Let
\begin{equation}
    \label{East}
    E^\ast(\tau, t)= \zeta^\ast(3t) E(\tau,t)
\end{equation}
where $\zeta^\ast(s)=\pi^{-s/2} \Gamma(s/2) \zeta(s)$
and $\zeta(s)$ is the Riemann zeta function. $E^\ast(\tau,t)$ as a
function of $t$ admits a meromorphic continuation to $\mathbb{C}$ with
poles at $t=1$, $0$ of residues $2/3$, $-2/3$ respectively.
It is shown in \cite{efrat2} that 
\begin{align}
  \label{East-expand}
  E^\ast(\tau, t) &= \frac{2/3}{t-1} + \gamma - \log 4\pi -
  \frac{2}{3} \log y_1 y_2^2 - 4 \log g(\tau) + O(t-1)
\end{align}
where $g(\tau)$ is given in \eqref{gtau}.

Let $Q$ be a 3 by 3 positive definite matrix associated with $\tau$ as in
\eqref{Qtau} so that
\[ Q=(\det Q)^{1/3} (y_1^2y_2)^{-2/3} \tau \tau'. \] 
One relates $Z_Q(s)$ to $E^\ast(\tau,t)$:
\begin{align}
  Z_Q(s) &= \sum_{m \in \mathbb{Z}^3 \backslash \{0\}} \frac{1}{Q[m]^s} \nonumber \\
  & = (\det Q)^{-s/3} \sum_{m \in \mathbb{Z}^3 \backslash \{0\}}
  \frac{(y_1^2y_2)^{2s/3}}{|m\tau|^{2s}} \nonumber \\
  & = (\det Q)^{-s/3} \sum_{m \in \mathbb{Z}^3 \backslash \{0\}}
  \frac{(y_1^2y_2)^{2s/3}}{ \left ( (m_1 y_1 y_2)^2 + (m_1 y_1 x_2 + m_2 y_1)^2
    + (m_1 x_3 + m_2 x_1 + m_3)^2 \right )^s} \nonumber \\
  &= (\det Q)^{-s/3} \sum_{m \in \mathbb{Z}^3 \backslash \{0\}}
  \frac{(y_1^2y_2)^{2s/3}}{ \left ( y_1^2  |m_1 z_2 + m_2|^2
    + (m_1 x_3 + m_2 x_1 + m_3)^2 \right )^s} \nonumber \\
  &= (\det Q)^{-s/3} \zeta(2s) \sum_{(m_1,m_2,m_3)=1} 
  \frac{(y_1^2y_2)^{2s/3}}{ \left ( y_1^2  |m_1 z_2 + m_2|^2
    + (m_1 x_3 + m_2 x_1 + m_3)^2 \right )^s} \nonumber \\
  & = (\det Q)^{-s/3} \zeta(2s) E \left (\tau, \frac{2s}{3} \right )
    \nonumber \\
& = (\det Q)^{-s/3} \pi^s \Gamma(s)^{-1}
  E^\ast \left (\tau, \frac{2s}{3} \right )\label{ZQtoEast}
\end{align}
By \eqref{East-expand} with $t=2s/3$, \eqref{ZQtoEast} becomes
\begin{align}
Z_Q(s) &= (\det Q)^{-s/3} \pi^s \Gamma(s)^{-1} \left( \frac{1}{s-3/2} + \gamma -
\log 4\pi -\frac{2}{3} \log y_1 y_2^2 - 4 \log g(\tau)
+ O\left(s-\frac{3}{2} \right ) 
  \right) \label{expand0}
\end{align}
Expand the other terms in \eqref{ZQtoEast} about $s=3/2$ to find
\begin{align}
  (\det Q)^{-s/3} \pi^s &= \pi^{3/2} (\det Q)^{-1/2} + \pi^{3/2} (\det Q)^{-1/2}
  \log \left ( \pi (\det Q)^{-1/3} \right ) \left( s-\frac{3}{2} \right)
  + O\left(s-\frac{3}{2} \right ) \label{expand1}
  \\ \Gamma(s) &= \Gamma \left( \frac{3}{2} \right )
  + \Gamma' \left ( \frac{3}{2} \right ) \left ( s - \frac{3}{2} \right )
  + O\left(s-\frac{3}{2} \right )  \label{expand2}
\end{align}
Note that 
\[ \Gamma \left( \frac{3}{2} \right ) = \frac{\pi^{1/2}}{2}, \ \ 
\frac{\Gamma'  (3/2) }{\Gamma (3/2)}= \psi\left(\frac{3}{2} \right)
= 2- 2\log 2 - \gamma \]
where $\psi = \Gamma'/\Gamma$ is the digamma function. Consequently
by \eqref{expand1}, and \eqref{expand2}, \eqref{expand0}
simplifies to
\begin{equation}
  Z_Q(s) =  \frac{2\pi (\det Q)^{-1/2}}{s-3/2} + 
2\pi (\det Q)^{-1/2}  \left ( \log(\det Q)^{-1/3} + 2 \gamma -2 -
       2 \log y_1^{1/3} y_2^{2/3}
       - 4 \log g(\tau)  \right ) + O\left( s-\frac{3}{2} \right )
\end{equation}
and the lemma follows.
\end{proof}

\begin{lemma}
  \label{l-R0k3}
  Let $\Lambda$ be a three-dimensional lattice. Then
  \[ k_3(Q^\ast) = 2\pi (\det Q^\ast)^{-1/2} \left ( 2\gamma -2 +2\log(2\pi)
  + (2\pi)^2 R_{\Lambda,3/2}(0) \right ) \] 
where $Q^\ast$ is a Gram  matrix associated with $\Lambda^\ast$. 
\end{lemma}

\begin{proof}
   Compare $R_{\Lambda.3/2}(0)$ in Lemma \ref{l-1}
   and $k_3(Q)$ in Lemma \ref{l-k3} with $Q$ substituted by $Q^\ast$.
  \end{proof}

Recall the height $h(M)$ of a Riemannian manifold defined in \eqref{height}. 
When $M$ is a flat torus $\Tn$, the eigenfunctions
of $-\Delta$ are $e^{2\pi ik\cdot w}$, $k \in \Lambda^\ast$, and the corresponding
eigenvalues are $4 \pi^2 |k|^2$. Then the height of $\Tn$ is essentially
given by Epstein's zeta function because 
\begin{equation}
  Z(\Tn,s) = \sum_{k \in \Lambda^\ast \backslash \{0\}} \frac{1}{(2\pi)^{2s} |k|^{2s}}
  = (2\pi)^{-2s} \sum_{m \in \mathbb{Z}^n \backslash \{0\}} \frac{1}{Q^\ast[m]^s}
  = (2\pi)^{-2s} Z_{Q^\ast}(s) \label{Ztorus}
\end{equation}
where $Q^\ast$ is a Gram matrix associated with the dual lattice $\Lambda^\ast$.

We reprove a result in \cite[Chiu]{chiu}, correcting a minor error. 

\begin{lemma}
\label{l-hk}
  Let $Q$ be a Gram matrix associated with lattice $\Lambda$. Then
  \[ h(\Tn) = (\det Q)^{1/2} \pi^{-d/2} \Gamma\left (\frac{d}{2} \right )
  k_d(Q) +\psi\left( \frac{d}{2} \right )
   -\gamma + 2\log 2 \]
   where $\psi = \Gamma' / \Gamma$ is the digamma function. In particular
   \[ h(\mathbb{R}^3/\Lambda^\ast) = (\det Q^\ast)^{1/2} (2\pi)^{-1} k_3(Q^\ast) - 2 \gamma +2
   \]
\end{lemma}

\begin{proof}
By \eqref{Ztorus}, the height of $\Tn$ is
\begin{align}
  h(\Tn)&= Z'(\Tn, 0) \nonumber \\
  &= -2\log (2\pi) Z_{Q^\ast}(0) + Z_{Q^\ast}'(0). \label{hZQast}
\end{align}
By Lemma \ref{l-functional}, 
\begin{equation}
  Z_{Q^\ast}(s) = (\det Q^\ast)^{-1/2} \pi^{2s-d/2} \Gamma(s)^{-1}
  \Gamma\left(\frac{d}{2}-s\right ) Z_Q\left (\frac{d}{2}-s \right). 
\end{equation}
For $s$ near $0$, expand
\begin{align}
 \pi^{2s-d/2} &= \pi^{-n/2} +  \left ( 2\pi^{-n/2} \log \pi \right) s + O(s^2) \\
   \Gamma(s)^{-1} &= s+ \gamma s^2 + O(s^3) \\
   \Gamma\left(\frac{d}{2}-s\right) &= \Gamma\left( \frac{d}{2} \right)
   - \Gamma'\left (\frac{d}{2} \right ) s + O(s^2) \\
   Z_Q\left( \frac{d}{2}-s \right ) &= -\frac{(\det Q)^{-1/2} \pi^{d/2}
      \Gamma\left (\frac{d}{2} \right )^{-1}}{s} + k_d(Q) + O(s) \label{line4}
\end{align}
where \eqref{line4} follows from \eqref{k}. Then 
\begin{align}
  Z_{Q^\ast}(s) &= -1 + \left ( (\det Q)^{1/2} \pi^{-d/2}
  \Gamma\left(\frac{d}{2} \right ) k_d(Q)
   + \psi\left (\frac{d}{2} \right) - \gamma -2\log \pi\right ) s + O(s^2) 
\end{align}
Hence 
\begin{align}
  Z_{Q^\ast}(0) &= -1 \\
  Z_{Q^\ast}'(0) &=  (\det Q)^{1/2} \pi^{-d/2}
  \Gamma\left(\frac{d}{2} \right ) k_d(Q)
   + \psi\left (\frac{d}{2} \right) - \gamma -2\log \pi
\end{align} 
and consequently by \eqref{hZQast}
\begin{align}
  h(\Tn) = (\det Q)^{1/2} \pi^{-d/2} \Gamma\left ( \frac{d}{2} \right ) k_d(Q)
                 + \psi \left ( \frac{d}{2} \right ) - \gamma +2\log 2
\end{align}
from which the lemma follows.
\end{proof}

We now state the key lemma in this paper.

\begin{lemma} \label{l-2}
  If $d=3$, then
  \[ h(\mathbb{R}^3\slash\Lambda^\ast) = (2\pi)^2 R_{\Lambda,3/2}(0)
  + 2\log(2\pi)  \]
\end{lemma} 

\begin{proof}
  Compare $R_{\Lambda,3/2}(0)$ in Lemma \ref{l-R0k3} and
  $h(\mathbb{R}^3/\Lambda^\ast)$ in Lemma \ref{l-hk}.
  \end{proof}

To prove the main theorem, we need a deep result of Sarnak and Str\"ombergsson.

\begin{theorem}[{\cite[Sarnak and Str\"ombergsson]{sarnak-strombergsson}}]
  Among three dimensional lattices $\Lambda$ of unit volume, the height
  $h(\T)$ is uniquely minimized by the FCC lattice. \label{t-SS}
  \end{theorem}

\begin{proof}[Proof of Theorem \ref{t-optimal}]

  By Lemma \ref{l-2} and Theorem \ref{t-SS}, $R_{\Lambda, 3/2}(0)$ is
  uniquely minimized when $\Lambda^\ast$ is the FCC lattice.
  Since the FCC lattice is the dual lattice of the BCC lattice, $\Lambda$
  must be the BCC lattice. 
\end{proof}

\bibliography{citation}
\bibliographystyle{plain}

\end{document}